\documentstyle[12pt]{article}
\newenvironment{proof}{\noindent {\em Proof. }}{\hfill$\Box$.\\}
\newtheorem{theo}{Theorem}
\newtheorem{prop}[theo]{Proposition}
\newtheorem{lemm}[theo]{Lemma}
\newtheorem{coro}[theo]{Corollary}
\newtheorem{defi}[theo]{Definition}
\newtheorem{rema}[theo]{Remark}
\newcommand{\e}{\label}
\newcommand{\r}[1]{(\ref{#1})}
\newcommand{\re}{\ref}
\newcommand{\k}{\ldots}
\newcommand{\p}{\partial}
\newcommand{\di}{\p\!\!\!/}
\newcommand{\Bo}{\mbox{\raisebox{-0.3ex}{\Large\mbox{$\Box$}}}}
\newcommand{\api}{{\cal A}}
\newcommand{\bpi}{{\cal B}}
\newcommand{\cpi}{{\cal C}}
\newcommand{\fpi}{{\cal F}}

\newcommand{\ppi}{{\cal P}}
\newcommand{\spi}{{\cal S}}

\newcommand{\gpi}{{\cal G}}
\newcommand{\vhi}{\varphi}

\newcommand{\rb}{\mbox{{\bf R}}{}}
\newcommand{\cb}{\mbox{{\bf C}}{}}
\newcommand{\jb}{\mbox{{\bf 1}}{}}

\newcommand{\ga}{\gamma}

\newcommand{\eps}{\varepsilon}

\newcommand{\lam}{\lambda}
\newcommand{\si}{\sigma}

\newcommand{\de}{\Delta}
\newcommand{\la}{\Lambda}
\newcommand{\wb}{\bar w}

\newcommand{\g}{\Gamma}

\newcommand{\qpg}{quantum Poincar\'e  group}

\newcommand{\qms}{quantum Minkowski space}
\newcommand{\be}{\begin{equation}}
\newcommand{\ee}{\end{equation}}
\newcommand{\bt}{\begin{theo}}
\newcommand{\et}{\end{theo}}
\newcommand{\bp}{\begin{prop}}
\newcommand{\ep}{\end{prop}}
\newcommand{\bl}{\begin{lemm}}
\newcommand{\el}{\end{lemm}}
\newcommand{\bc}{\begin{coro}}
\newcommand{\ec}{\end{coro}}
\newcommand{\bde}{\begin{defi}}
\newcommand{\ede}{\end{defi}}
\newcommand{\br}{\begin{rema}}
\newcommand{\er}{\end{rema}}
\newcommand{\bd}{\begin{proof}}
\newcommand{\ed}{\end{proof}}
\newcommand{\ba}{\begin{array}}
\newcommand{\ea}{\end{array}}
\newcommand{\btr}{\begin{trivlist}}
\newcommand{\etr}{\end{trivlist}}
\newcommand{\lra}{\longrightarrow}
\newcommand{\ot}{\otimes}
\newcommand{\op}{\oplus}
\newcommand{\ti}{\times}
\newcommand{\sd}{\rhd\mbox{\hspace{-2ex}}<}              
\newcommand{\spa}{\mbox{{\rm{span}}}}                 
                 
\newcommand{\po}{\mbox{{\rm{Poly}}}}
\newcommand{\id}{\mbox{{\rm{id}}}}

\newcommand{\te}{\tilde}

\newcommand{\Rep}{\mbox{{\rm{Rep}}}\ }

\newcommand{\ov}{\overline}

\newcommand{\tp}{\otimes}
\newcommand{\poi}{Poincar\'e\ }
\newcommand{\w}{\wedge}
\begin{document}
\title{Quantum Minkowski spaces}
\author{P. Podle\'s\\
Department of Mathematical Methods in Physics\\ 
Faculty of Physics, Warsaw University}
\date{}
\maketitle

One of the main problems of theoretical physics is to find a
satisfactory theory which would generalize both the quantum
field theory and the general theory of relativity. It is widely
recognized that in such a theory the geometry of the space--time
should drastically change at very small distances, comparable
with the Planck's length. One of possibilities is to replace
the space--time by so called quantum space. In such an approach
the role of the commutative algebra of functions on the
space--time (generated in the simplest case by the coordinates)
is played by some noncommutative algebra. 
At the present stage we are only able to test this idea in
particular examples, which are important in physics. Many
classical objects were already deformed in the above sense (cf.
e.g. \cite{W1}, \cite{P1}, \cite{FRT}, \cite{PW}, \cite{CSSW}, \cite{CSW},
\cite{WZ}). Here we deal with the most
interesting case, namely that of Minkowski space $M$ endowed
with the action of Poincar\'e group $P$. Examples of quantum
Poincar\'e groups and their actions on quantum Minkowski spaces
appeared e.g. in \cite{L}, \cite{ChD}, \cite{GQ}
(the case of quantum Poincar\'e groups and algebras 
extended by  dilatations was considered e.g. in \cite{D},
\cite{S}, \cite{OSWZ}, \cite{M}), see also \cite{Cas}, \cite{Az}. 
The aim is  to find the
classification of \qpg s and \qms s as well as mathematical and
physical properties of those objects. We sketch the results 
of four papers \cite{POI}, \cite{KG}, \cite{QUA}, \cite{GAM}.

The (connected component of) vectorial \poi group
\[ \te P=SO_0(1,3)\sd\rb^4=\{(M,a): M\in SO_0(1,3), a\in\rb^4\} \]
has the multiplication $(M,a)\cdot(M',a')=(MM',a+Ma')$. By the \poi group we
mean spinorial \poi group (which is more important in quantum field theory
then $\te P$)
\[ P=SL(2,\cb)\sd\rb^4=\{(g,a): g\in SL(2,\cb), a\in\rb^4\} \]
with multiplication $(g,a)\cdot(g',a')=(gg',a+\lam_g(a'))$ where the double
covering $SL(2,\cb)\ni g\lra \lam_g\in SO_0(1,3)$ is the
standard one.
The group
homomorphism $\pi:P\ni(g,a)\lra(\lam_g,a)\in\te P$ is also a double
covering. In particular, $(-\jb_2,0)\in P$ can be treated as rotation
about $2\pi$ which is trivial in $\te P$ but nontrivial in $P$ (it changes
the sign of wave functions for fermions). Both $P$ and $\te P$ act on
Minkowski space $M=\rb^4$ as follows $(g,a)x=(\lam_g,a)x=\lam_gx+a$, $g\in
SL(2,\cb)$, $a,x\in\rb^4$.

Let us consider continuous functions $w_{AB}$,
$y_i$ on $P$ defined by 
\[ w_{AB}(g,a)=g_{AB},\qquad y_i(g,a)=a_i. \]
We introduce 
Hopf ${}^*$-algebra $\po(P)=(\bpi,\de)$ of polynomials on the \poi
group $P$ as the 
${}^*$-algebra $\bpi$ with identity $I$, generated by $w_{AB}$
and $y_i$, $A,B=1,2$, $i\in\spi=\{0,1,2,3\}$ endowed with the
comultiplication $\de$ given by $(\de f)(x,y)=f(x\cdot y)$, $f\in\bpi$,
$x,y\in P$ ($f^*(x)=\ov{f(x)}$). In particular, 
\be \de w_{CD}=w_{CF}\ot w_{FD}, \e{1.1}\ee 
\be \de y_i=y_i\ot I+\la_{ij}\ot y_j, \e{1.1'}\ee
$y_i^*=y_i$, where 
\be \la=V^{-1}(w\tp\wb)V, \qquad V=\left(\ba{cccc} 
1&0&0&1 \\ 0&1&-i&0 \\ 0&1&i&0 \\ 1&0&0&-1 \ea \right) \e{1.2}\ee
(we sum over repeated indices, one has $V_{CD,i}=(\si_i)_{CD}$ where 
$\si_i$ are the Pauli matrices).
Moreover, $w$ is a representation, i.e. it is invertible (as
$2\times 2$ matrix with elements in $\bpi$) and satisfies \r{1.1}.
Equivalences of representations are defined as usual (by means
of invertible matrices with complex entries), e.g. $\la\simeq w\ot\wb$.
We put $y=(y_i)_{i\in\spi}$. One can also treat $w_{CD}$ as
continuous functions on the Lorentz group $L=SL(2,\cb)$ ($w_{CD}(g)=g_{CD}$,
$g\in L$). We define 
Hopf ${}^*$-algebra $\po(L)=(\api,\de)$ of polynomials on
$L$ as ${}^*$-algebra with $I$, generated by all $w_{CD}$ endowed with $\de$
obtained by restriction of $\de$ for $\bpi$ to $\api$. Clearly $w$ and $\la$
are representations of $L$. It is easy to check that
\begin{trivlist}
\item[1.]  $\bpi$ is generated as algebra by $\api$
and the elements $y_i$, $i\in\spi$. 
\item[2.] $\api$ is a Hopf ${}^*$-subalgebra of $\bpi$.
\item[3.] $\ppi=\left(\ba{cc}\la&y\\0&I\ea\right)$ is a representation of
$\bpi$ where $\la$ is given by \r{1.2}.
\item[4.] There exists $i\in\spi$ such that $y_i\not\in\api$.
\item[5.] $\g\api\subset\g$ where $\g=\api X+\api$, $X=\spa\{y_i:i\in\spi\}$.
\item[6.] The left $\api$-module $\api\cdot\spa\{y_iy_j, y_i, I:
i,j\in\spi \}$ has a free basis consisting of $10+4+1$ elements.
\end{trivlist}
( 5. and 6. follow from the relations $y_ia=ay_i$, $y_iy_j=y_jy_i$,
$a\in\api$, and elementary computations, a free basis is given by $\{y_iy_j,
y_i, I: i\leq j, i,j\in\spi \}$). According to \cite{WZ}, $\po(L)$
satisfies: 
\begin{trivlist}
\item[i.] $(\api,\de)$ is a Hopf ${}^*$-algebra such that $\api$ is generated 
(as
${}^*$-algebra) by matrix elements of a two--dimensional representation $w$
\item[ii.] $w\tp w\simeq I\oplus w^1$ where $w^1$ is a representation
\item[iii.] the representation $w\tp\wb\simeq\wb\tp w$ is irreducible
\item[iv.] if $\api',\de',w'$ satisfy i.--iii. and there exists a 
Hopf ${}^*$-algebra 
epimorphism $\rho:\api'\lra\api$ such that $\rho(w')=w$ then $\rho$ 
is
an isomorphism (the universality condition).
\end{trivlist}

We define quantum Lorentz groups and quantum \poi
groups 
as objects having the same properties as the classical
Lorentz and \poi groups:

\bde\e{d1.1} 
We say \cite{WZ} that $H$ is a quantum Lorentz group if $\po(H)=(\api,\de)$
satisfies i.--iv.
We say \cite{POI} that $G$ is a \qpg\ if 
Hopf ${}^*$-algebra $\po(G)=(\bpi,\de)$ satisfies the
conditions 1.--6. for some quantum Lorentz group $H$ with
$\po(H)=(\api,\de)$ and a representation $w$ of $H$. \ede

{\bf Remark.} The condition 5. follows from $\ppi\tp w\simeq w\tp\ppi$,
$\ppi\tp\wb\simeq\wb\tp\ppi$, while 6. is suggested by the requirement
$W(\ppi\tp\ppi)=(\ppi\tp\ppi)W$ for a ``twist-like'' matrix $W$. Moreover, 
the condition 4. is superfluous (it follows from the 
condition 6.).

\bt\e{t1.2} Let $G$ be a \qpg, $\po(G)=(\bpi,\de)$. Then $\api$ is 
linearly 
generated by matrix elements of irreducible representations of $G$, 
so
$\api$ is uniquely determined. Moreover, we can choose $w$ in such a way that
$\api$ is the universal
${}^*$-algebra generated by $w_{AB}$, $A,B=1,2,$ satisfying
\[ (w\tp w)E=E, \]
\[ E'(w\tp w)=E', \]
\[ X(w\tp\wb)=(\wb\tp w)X, \]
where the triples $E\in M_{4\times 1}(\cb),
E'\in M_{1\times 4}(\cb), 
X\in M_{4\times 4}(\cb)$ are listed in Theorem 1.4 of \cite{POI}.
We can (and will) choose $y_i$ in such a way that $y_i^*=y_i$.\et   

In particular, it turns out that only quantum Lorentz groups
with the parameter $q=\pm1$ are admissible.
Nevertheless, there
are many families of admissible quantum Lorentz groups (numbered
by some other parameters).
In the following we assume that $G$ is a \qpg, $\po(G)=(\bpi,\de)$ and $w,y$
are as in Theorem \re{t1.2}. Using
the general theory of inhomogeneous quantum groups \cite{INH},
we find the full system of commutation relations for
$\bpi$:

\bt\e{t1.3} $\bpi$ 
is the universal ${}^*$-algebra with $I$, generated by $\api$ and $y_i$
with relations
\[ y\ot v=G_v(v\ot y)+H_vv-(\la\ot v)H_v,\]
\[ (R-\jb)(y\ot y-Zy+T-(\la\ot\la)T)=0,\]
\[ y_i^*=y_i,\quad i\in\spi,\]
for any $v\in\Rep G$ (the set of all representations of $G$), where
$(G_v)_{iC,Dj} = f_{ij}(v_{CD})$, $(H_v)_{iC,D} = \eta_i(v_{CD})$,
$R=G_{\la}$, $Z=H_{\la}$ and 
$T_{ij}\in\cb$, $f_{ij},\eta_i\in\api'$ ($i,j\in\spi$) satisfy
the conditions of Theorem 1.5 of \cite{POI}. 
$\de$ is given by \r{1.1},\r{1.1'}.
Moreover, $(\bpi,\de)$ gives a \qpg\ if and only if a system of
linear and quadratic equations listed in the proof of Theorem 1.6 of \cite{POI}
is fulfilled.\et

It turns out that there are two choices for $f$
determined by a number $s=\pm1$ - the calculations are made for
each $s$ separately and the results are given in terms of 
$H_{EFCD}=V_{EF,i}\eta_i(w_{CD})$ and $T_{EFCD}=V_{EF,i}V_{CD,j}T_{ij}$.
We solve \cite{POI} the above  system of equations for
almost all quantum Lorentz groups (except two cases, including
the classical Lorentz group for which a large class of solutions
is known). 
Moreover, we single out unisomorphic objects.
The classification is presented in Theorem 1.6 of
\cite{POI}. 
We also identify few
examples which were known earlier (cf. \cite{L}, \cite{ChD}, \cite{GQ}).
We prove that $\bpi$ 
has exactly the sime ``size'' as in the undeformed case.
Namely,
\[ \bpi^N=\api\cdot\spa\{y_{i_1}\cdot\k\cdot y_{i_n}: i_1,\k,i_n\in\spi,\quad
n=0,1,\k,N \} \]
is a free left $\api$-module and
\[ dim_{\api}\bpi^N=\sum_{n=0}^N d_{n} \]
where $d_n$ is the number of classical 
monomials of $n$th degree in $4$ variables.

We denote by $l:P\times M\lra M$ the action of \poi group on Minkowski space,
$\cpi=\po(M)$ denotes the unital algebra generated by coordinates $x_i$
($i\in\spi$) of the Minkowski space $M=\rb^4$. 
The coaction $\Psi:\cpi\lra\bpi\ot\cpi$ and $*$ in $\cpi$ 
are
given by $(\Psi f)(x,y)=f(l(x,y))$, 
$f^*(y)=\ov{f(y)}$, $x\in P$, $y\in M$. One gets
\be \Psi x_i=\la_{ij}\ot x_j+y_i\ot I.\e{1.16}\ee

We define a quantum Minkowski space as object having the same
properties as the classical Minkowski space:

\bde\e{d1.6} We say that $(\cpi,\Psi)$ describes a quantum Minkowski space
associated with a \qpg\ $G$, $\po(G)=(\bpi,\de)$, if
$\cpi$ is a unital 
${}^*$-algebra generated by $x_i$, $i\in\spi$,
$\Psi:\cpi\lra \bpi\ot\cpi$ is a unital ${}^*$-homomorphism, \r{1.16}
holds and other conditions of Definition 1.10 of \cite{POI} are satisfied.\ede

\bt\e{t1.7} Let $G$ be a \qpg\ with $w,y$ as in Theorem 
\re{t1.2}. 
Then there exists a unique (up to a ${}^*$-isomorphism) pair $(\cpi,\Psi)$
 describing
associated Minkowski space:

$\cpi$ is the 
universal unital ${}^*$-algebra generated by $x_i$, $i=0,1,2,3$, satisfying
${x_i}^*=x_i$ and 
\[ (R-\jb)(x\ot x-Zx+T)=0,\]
and $\Psi$ is given by $\r{1.16}$. Moreover, 
\[ \dim\cpi^N=\sum^N_{n=0} d_n, \]
where $\cpi^N=\spa\{x_{i_1}\cdot\k\cdot x_{i_n}:\ \ i_1,\k,i_n\in\spi,\ \ 
n=0,1,\k,N\}.$\et

The next our goal is to find the differential structure on quantum
Minkowski spaces. It turns out \cite{KG} that there exists a
unique $4$-dimensional  covariant first order differential
calculus on a quantum
Minkowski space provided  $\te F=0$ where
\[ \te F= [(R-\jb)\ot\jb]\{(\jb\ot Z)Z-(Z\ot\jb)Z+T\ot\jb-
 (\jb\ot R)(R\ot\jb)(\jb\ot
T)\} \]
(otherwise there is no such a calculus).
This condition singles out a large class of quantum Minkowski spaces
\cite{KG}. From now on we assume that this condition is
fulfilled. 
 In particular, there exists a
$\cpi$-bimodule $\Gamma^{\w 1}$ (of differential forms of the first order)
and a linear mapping $d: {\cal C} \lra \Gamma^{\w 1}$ such that
$d(ab) = a(db) + (da)b$, $a,b \in {\cal C}$, and
$dx_i$, $i\in\spi$,
form a basis of $\Gamma^{\w 1}$ (as right ${\cal C}$-module). 
We prove
\[ x_idx_j = R_{ij,kl} dx_kx_l + Z_{ij,k}dx_k,\
i,j\in\spi. \]
 This calculus prolongates to a unique exterior algebra of differential
forms, with the same properties as in the undeformed case. In particular
it possesses $*$ such that $(dx_i)^*=dx_i$.

The
partial derivatives 
$\p_i:
{\cal C} \lra {\cal C}$
are uniquely defined by
\[ da = dx_i\p_i(a), \quad a\in\spi. \]
They can be also obtained as 
\[ \p_i = (Y_i \otimes id)\Psi \]
where $Y_i\in\api'$ are introduced
in the proof of Proposition 3.1.2 of \cite{KG}.

The metric tensor $g = (g_{ij})_{i,j\in\spi}$ (its entries are
called in \cite{KG} by $g^{ij}$) is defined by the equations
$(\la\ot\la)g=g$ (or $\la g\la^T=g$) and $\ov{g_{ij}}=g_{ji}$. 
After the choice of a real factor we fix it as
\[ g=-2q^{1/2}(V^{-1}\ot V^{-1})(\jb\ot X\ot\jb)(E\ot\tau E),
\]
where $\tau$ is always a standard twist. Then the Laplacian is
defined by $\Bo=g_{ij}\p_j\p_i$. It commutes with the partial
derivatives and therefore the momenta $P_l=i\p_l$ are well
defined in the spaces of solutions of the Klein--Gordon equation
$(\Bo+m^2)\vhi=0$. One proves  that the
momenta\footnote{in this case we relax our rule of writing all
indices in the subscript position}
$P^k=g_{kl}P_l$ and
Laplacian are hermitian and have good transformation properties.
The commutation relations among partial derivatives and with the
coordinates are as follows:
\[ \p_l\p_k = R_{ij,kl}\p_j\p_i,\]
\[ \p_ix_k = \delta_{ki} + (R_{kl,in} x_n + Z_{kl,i})\p_l. \]

Let us now pass to the particles of spin $1/2$ \cite{GAM}.
First we define the space of bispinors as $\cb^4$
endowed with a representation $\gpi\simeq w\oplus\wb$ of a \qpg.
We choose $\gpi={}^cw\op\wb$ 
where ${}^cw=(w^T)^{-1}\simeq w$. 
We are going to find the gamma matrices
$\ga_i\in M_{4\ti 4}(\cb)$,
$i\in\spi$. At the moment they are not determined yet. The
Dirac operator has form $\di=\ga_i\ot\p_i$. It acts on
the bispinor functions 
$\phi\in\te\cpi\equiv\cb^4\ot\cpi$ (in a more
advanced approach we should consider square integrable functions
$\phi$). They can be written as
$\phi=\eps_a\ot\phi_a$ where $\eps_a$, $a=1,2,3,4$, form
the standard basis of $\cb^4$. We define the action
$\te\Psi:\te\cpi\lra\bpi\ot\te\cpi$ of \qpg\ on $\te\cpi$ by 
\[
\te\Psi(\eps_a\ot\phi^a)=\gpi_{al}\phi_a^{(1)}\ot\eps_l\ot\phi_a^{(2)}
\]
where $\Psi(\phi^a)=\phi_a^{(1)}\ot\phi_a^{(2)}$
(Sweedler's notation, exception of the summation
convention). Then the classical condition of invariance of the
Dirac operator is generalized to 
\[ \te\Psi(\di\phi)=(\id\ot\di)[\te\Psi(\phi)],\quad \phi\in\te\cpi.
\]
According to Theorem III.1 of \cite{GAM}, the above condition is 
equivalent to
\[ \ga_i=\left(\ba{cc} 0&bA_i\\
                a\si_i&0\ea\right),\quad i\in\spi, \]
where
\be A_i=q^{-1/2}E^T(\si_i\circ D)E,\e{2.12}\ee
$(\si_i\circ D)_{KL}=(\si_i)_{AB}D_{AB,KL}$, $D=\tau X^{-1}\tau$,
$a,b\in\cb$
($E$ is regarded here as $2\ti 2$ matrix). Thus we have found
the form of the Dirac operator up to two constants. But Theorem
III.2 of \cite{GAM} says that the following are equivalent:
\btr
\item[1.] $\di^2=\Bo$.
\item[2.] $\ga_i\ga_j+R_{ji,lk}\ga_k\ga_l=2g_{ji}\jb,\ \ \ \ \
i,j\in\spi$. 
\item[3.] $ab=1$.
\etr
Moreover, the remaining freedom in the choice of  $a$ results in a
trivial scaling of the undotted spinor and we can set $a=b=1$.
Thus we have obtained
\[ \ga_i=\left( \ba{cc} 0&A_i\\ \si_i&0\ea\right)  \]
where $A_i$ are given by \r{2.12}. Now the form of the Dirac equation
$(i\di+m)\vhi=0$ is determined.

In the next step we 
find (formal) solutions of Klein--Gordon and Dirac equations in
two important cases (Section 4 of \cite{KG}). In one of them we use (as a
tool during calculations) an additional algebra $\fpi$ and its
representations. Specific calculations are made in \cite{GAM}
(including the form of metric tensor, gamma matrices,
representations of $\fpi$ and the solutions of Klein--Gordon and
Dirac equations). For spin $0$ particles the momenta $P_j=i\p_j$. 
For spin $1/2$
particles we set \cite{GAM} the momenta as
\be \te P_j=i\te\p_j, \e{3.7}\ee
\be \te\p_j=(Y_j\ot\id)\te\Psi:\te\cpi\lra\te\cpi \e{3.8}\ee
(motivation: \r{3.7}--\r{3.8} remain true if we omit tildas everywhere).
We get four objects $\te P^k=g_{kj}\te P_j$ 
which have good transformation properties, 
commute with the Dirac
operator and are hermitian w.r.t. (indefinite) inner product in
the space of bispinors. However, in many cases their spectral
properties are not satisfactory. The problem of further improvement in
this matter remains open. We also study certain expressions like
the deformed Lagrangian. They transform themselves in a similar
way as in the standard theory.

It turns out (cf. Theorem 1.13 of \cite{POI}) that there exist
invertible matrices $W$ such that $W(\ppi\ot\ppi)=(\ppi\ot\ppi)W$ and the
Yang--Baxter equation
\[ (W\ot\jb)(\jb\ot W)(W\ot\jb)=(\jb\ot W)(W\ot\jb)(\jb\ot W) \]
is satisfied. Namely, up to a constant they are given by the
unit matrix and
\[ R_Q=\left(\ba{cccc} R & Z & -R\cdot Z & (R-\jb^{\ot 2})T+b\cdot g \\
                        0 & 0 &      \jb  &     0      \\
                        0 &\jb&      0    &     0      \\
                        0 & 0 &      0    &     1      \ea\right),\]
where $b\in\cb$.

In \cite{QUA} we
show the existence of 
${\cal R} \in ({\cal B} \otimes {\cal B})'$
such that $({\cal B},\de,{\cal R})$ is a coquasitriangular (CQT) 
Hopf algebra. In other words,  for any $v,z\in\Rep G$ 
we define
\[ R^{vz} \in
\mbox{ Lin}({\cb}^{\dim v} \otimes {\cb}^{\dim z},{\cb}^{\dim z}
\otimes {\cb}^{\dim v})\]
 by
\[ (R^{vz})_{ij,kl} = {\cal R}(v_{jk} \otimes z_{il}),\ 
j,k = 1,\dots,\dim v,\ i,l
= 1,\dots,\dim z, \]
and prove
\[ R^{1v} = R^{v1} = \jb,\]
\[ R^{v_1 \otimes v_2,z} = (R^{v_1z} \otimes \jb)(\jb \otimes R^{v_2z}),\]
\[ R^{v,z_1\otimes z_2} = (\jb \otimes R^{vz_2})(R^{vz_1} \otimes \jb),\]
\[ (z\ot v)R^{vz}=R^{vz}(v\ot z),\]
for all $v,v_1,v_2,z,z_1,z_2\in\Rep G$ ($1=(I)$ is the trivial representation).

The classification of all CQT Hopf algebra structures (for
all quantum Poincar\'e groups) is done in Theorem 3 of \cite{QUA}.
In particular,  $R^{\ppi\ppi}=R_Q$,
\[ R^{v\ppi} = \left( \begin{array}{c} G_v,\ H_v \\
0,\ \jb
\end{array} \right),\]
$R^{\ppi v}=(R^{v\ppi})^{-1}$, 
$R^{ww} = kL$, $R^{w{\bar w}} = kX$, $R^{{\bar w}w} = qkX^{-1}$, 
$R^{{\bar w}{\bar
w}} = k\tau L\tau$, for all representations $v$ of the
corresponding quantum Lorentz group $H$ (these data determine
${\cal R}$ uniquely), where
$L = sq^{1/2}(\jb+qEE')$, 
$k = \pm 1$ (two possible ${\cal R}$ for each $b \in {\cb}$).

We have to do with CQT Hopf ${}^*$-algebra iff 
 $\overline{{\cal R}(y^* \otimes x^*)} = {\cal R}(x \otimes y)$, $x,y
\in {\cal B}$, iff $q=1$ and $b\in\rb$. We have cotriangular (CT) Hopf
algebra iff $(R^{vz})^{-1} = R^{zv}$ for all $v,z \in\Rep G$ iff
$q=1$ and $b=0$ (so then it is also a CT Hopf ${}^*$--algebra).

Using the above results, universal enveloping algebras for quantum Poincar\'e
groups are introduced. Their commutation relations are
investigated. Moreover, we classify C(Q)T Hopf (${}^*$--)algebra 
structures for quantum Lorentz groups.
We also show
some general statements concerning coquasitriangularity. 
The results of \cite{QUA} are used in Section 5 of \cite{KG} to define the
Fock space for non-interacting particles of spin $0$ on a
quantum Minkowski space. 
Then we take a $CT$ Hopf ${}^*$-algebra structure
${\cal R}$ as 
above and 
introduce the particles interchange operator $K: {\cal C} \otimes {\cal C} \lra
{\cal C} \otimes {\cal C}$ by 
\[
K(x \otimes y) = {\cal
R}(y^{(1)} \otimes x^{(1)})(y^{(2)} \otimes x^{(2)}),
\] 
where 
$\Psi(x) = x^{(1)} \otimes x^{(2)}$, $\Psi(y) =
y^{(1)} \otimes y^{(2)}$. 
It defines the action of the permutation group $\Pi_n$ in
$\cpi^{\ot n}$ 
($\Pi_n\ni\sigma\lra\pi_{\sigma}$)
which agrees with the action of the \qpg\ G. Thus
$G$ acts in the boson subspace $\cpi^{\ot_s n}$. 
If $W: {\cal C} \to {\cal C}$ is an operator related to a
single particle then the corresponding $n$-particle operator is given by
\[
W^{(n)} = \sum_{m=1}^n \pi_{(1,m)} ( W\otimes \jb^{\otimes
(n-1)})\pi_{(1,m)}:
{\cal C}^{\otimes_s n} \to {\cal C}^{\otimes_s n}
\]
(the $m$-th term is the operator in ${\cal C}^{\otimes n}$ 
corresponding to the $m$-th particle).
We can also define the Fock space $F = \oplus_{n=0}^{\infty} {\cal
C}^{\otimes_s n}$ and the operator $\oplus_{n=0}^{\infty} W^{(n)}$ acting in
$F$.

For particles of mass $m$ we should consider $\ker(\Bo + m^2)$ instead of
${\cal C}$ and a scalar product there (heuristically e.g. $W=P^k$, $k\in\spi$, 
would be hermitian operators in such a space).

Results of \cite{KG} and \cite{QUA} are
proven also for general inhomogeneous quantum groups
(satisfying certain conditions). References to the existing
literature are given in \cite{POI}, \cite{KG}, \cite{QUA}, \cite{GAM}.

Concluding, quantum Minkowski spaces and quantum Poincar\'e
groups have a lot of properties similar to that of  the classical ones. It
suggests a possibility of building more advanced physical models
which use those objects. However, it would need further studies
concerning deformed quantum field theory and interaction of particles.
There is also another advantage of quantum Minkowski spaces:
the fact that there are many possibilities in the choice of
parameters somehow forces us 
to find the proofs which have good geometric
meaning. In particular, the invariance of the Dirac
operator turns out to be equivalent to the fact that some object built from
the gamma matrices intertwines two specific representations of
the quantum Poincar\'e group. Then the form of gamma matrices 
is easy to find (without using the Lie algebra at all).


\newpage

\begin{thebibliography}{xx}
\bibitem{Cas} Aschieri, P. and Castellani, L., R-matrix formulation
of the quantum inhomogeneous groups $ISO_{qr}(N)$ and $ISp_{qr}(N)$,
{\em Lett. Math. Phys.} {\bf 36} (1996), 197--211; Bicovariant calculus
on twisted $ISO(N)$, quantum Poincar\'e group and quantum Minkowski space,
{\em Int. J. Mod. Phys. A} {\bf 11} (1996), 4513; Aschieri, P., 
Castellani, L. and Scarfone, A.M., Quantum orthogonal planes: $ISO_{qr}(N)$
and $SO_{qr}(N)$ -- bicovariant calculi and differential geometry on
quantum Minkowski space, q-alg/9709032. 

\bibitem{Az} de Azc\'arraga, J.A. and Rodenas, F., Differential calculus on
$q$-Minkowski space, {\em An. Fisica (Monogr.)} {\bf 2} (1995), 107--130;
de Azc\'arraga, J.A., Kulish, P.P. and Rodenas, F., On the physical contents
of $q$-deformed Minkowski spaces, {\em Phys. Lett. B} {\bf 351} (1995), 123;
Twisted $h$-spacetimes and invariant equations, {\em Z. Phys. C} {\bf 76}
(1997), 567--576.

\bibitem{CSSW} Carow--Watamura, U., Schlieker, M., Scholl, M. and
Watamura, S., Tensor representation of the quantum group $SL_q(2,\cb)$ 
and quantum Minkowski space, {\em Z. Phys. C -- Particles and Fields}
{\bf 48} (1990), 159--165.

\bibitem{CSW} Carow--Watamura, U., Schlieker, M. and Watamura, S.,
$SO_q(N)$ covariant differential calculus on quantum space and quantum 
deformation of Schr\"odinger equation, {\em Z. Phys. C -- Particles and 
Fields} {\bf 49} (1991), 439--446.

\bibitem{ChD} Chaichian, M. and Demichev, A.P., 
Quantum Poincar\'e group,
{\it Phys. Lett.\/}
{\bf B304} (1993), 220--224. 
\bibitem{D} Dobrev, V.K., Canonical $q$-deformations of noncompact Lie (super-)
algebras, {\it J.Phys.A: Math. Gen.\/} {\bf 26}(1993), 1317--1334.
\bibitem{L} Lukierski, J., Nowicki, A. and Ruegg, H., New quantum Poincar\'e 
algebra and $\kappa$--deformed field theory, {\it Phys. Lett.\/}
{\bf B293} (1992), 344--352;
Zakrzewski, S., Quantum Poincar\'e group related to the $\kappa$-Poincar\'e
algebra, {\it J. Phys. A: Math. Gen.\/} {\bf 27} (1994), 2075--2082.
\bibitem{M} Majid, S., Braided momentum in the $q$-Poincar\'e group,
{\it J. Math. Phys.\/} {\bf 34} (1993), 2045--2058.
\bibitem{OSWZ} Ogievetsky, O., Schmidke, W.B., Wess, J. and Zumino, B.,
$q$-Deformed Poincar\'e algebra, {\it Commun. Math. Phys.\/} {\bf 150}
(1992), 495--518.
\bibitem{P1} Podle\'s, P., Quantum spheres, {\em Lett. Math.
Phys.} {\bf 14} (1987), 193--202;
The classification of differential structures
on quantum $2$-spheres, {\em Commun. Math. Phys.} {\bf 150} (1992), 167--179;
Quantization enforces interaction.  Quantum
mechanics of two particles on a quantum sphere, {\em Int. J. Mod. Phys. A},
{\bf 7}, Suppl. 1B (1992), 805--812.
\bibitem{PW} Podle\'s, P. and Woronowicz, S.L., Quantum 
deformation of Lorentz group, {\it Commun. Math. Phys.\/} {\bf 130} (1990),
381--431.
\bibitem{INH} Podle\'s, P. and Woronowicz, S.L., On the 
structure of inhomogeneous quantum groups, 
hep-th/9412058, {\it Commun. Math. Phys.\/} {\bf 185} (1997), 325--358.
\bibitem{POI} Podle\'s, P. and Woronowicz, S.L., On the
classification of quantum Poincar\'e groups, {\it Commun. Math. Phys.\/}
{\bf 178} (1996), 61--82.
\bibitem{KG} Podle\'s, P., Solutions of
Klein--Gordon and Dirac equations on quantum Minkowski spaces,
{\it Commun. Math. Phys.\/} {\bf 181} (1996), 569--585.
\bibitem{QUA} Podle\'s, P.,
Quasitriangularity and enveloping algebras for inhomogeneous
quantum groups, {\it J. Math. Phys.\/} {\bf 37} (1996), 4724--4737.
\bibitem{GAM} Podle\'s, P., The Dirac operator and gamma
matrices for quantum Minkowski spaces, {\it J. Math. Phys.\/}
{\bf 38} (1997), 4474--4491.
\bibitem{FRT} Reshetikhin, N. Yu., Takhtadzyan, L. A. and Faddeev, L. D.,
Quantization of Lie groups and Lie algebras, {\em Leningrad Math. J.} {\bf
1:1} (1990), 193--225.  Russian original:  {\em Algebra i analiz} {\bf 1:1}
(1989), 178--206.
\bibitem{S} Schlieker, M., Weich, W. and Weixler, R., Inhomogeneous
quantum groups, {\it Z. Phys. C. -- Particles and Fields\/} {\bf 53} (1992),
79-82.
\bibitem{W1} Woronowicz, S.L., Twisted $SU(2)$ group. An example
of a non--commutative differential calculus, {\em Publ. RIMS,
Kyoto University} {\bf 23} (1987), 117--181; 
Compact matrix pseudogroups, {\em
Commun. Math. Phys.} {\bf 111} (1987), 613--665.
\bibitem{WZ} Woronowicz, S.L. and Zakrzewski, S., Quantum deformations of the
Lorentz group. 
The Hopf ${}^*$-algebra level, {\it Comp. Math.\/} {\bf 90} (1994),
211--243.
\bibitem{GQ} Zakrzewski, S., Geometric quantization of Poisson groups -- 
diagonal and soft deformations, {\it Contemp. Math.\/} {\bf 179} (1994),
271--285.

\end{thebibliography}
\end{document}